\documentclass[12pt, letterpaper]{article}
\usepackage[margin=2.5cm]{geometry} 

\usepackage{amssymb}
\usepackage{amsmath}
\usepackage{graphicx}
\usepackage{physics}

\usepackage{cite}
\usepackage{url}

\newcommand{\at}[2][]{#1\bigg|_{#2}}

\usepackage{xcolor}
\newcommand{\change}[1]{{#1}}
\newcommand{\subs}[1]{\textbf{#1}}

\title{Entropic dynamics yields reciprocal relations}

\author{ Pedro Pessoa$^1$ \\
$^1$Department of Physics, University at Albany - SUNY \\  Albany, NY - USA}
\date{}

\begin{document}

\maketitle
\abstract{
Entropic dynamics is a framework for defining dynamical systems that is aligned with the principles of information theory.
In an entropic dynamics model for motion on a statistical manifold,  we find that the rate of changes for expected values is linear to the gradient of entropy with reciprocal (symmetric) coefficients.
Reciprocity principles have been useful in physics since Onsager.
Here we show how the entropic dynamics reciprocity is a consequence of the information geometric structure of the exponential family, hence it is a general property that can be extended to a broader class of dynamical models.
}


\newpage

\section{Introduction}
Entropic dynamics (EntDyn) is a method to characterize dynamical processes that is aligned to the principles of information theory  taking into account the differential geometric structure of probability distributions given by information geometry (IG) \cite{Nielsen11,Caticha15,Amari16,Nielsen20}.
EntDyn has been successful in deriving dynamical processes in fundamental physics \cite{Caticha10,Caticha18,Ipek19}, but has also been expanded to other fields such as renormalization groups  \cite{Pessoa18}, finance \cite{Abedi19a}, machine learning \cite{NCaticha20} and network science \cite{Costa21} \change{some advances have also been made towards EntDyn applications in ecology \cite{Pessoa21}}.
Particularly in \cite{Pessoa20} EntDyn is presented as a general method for dynamical systems of probability distributions in the exponential family -- also referred to as Gibbs or canonical distributions.
Under a particular choice of sufficient statistics, the space of exponential families is parametrized by the expected values of these sufficient statistics  and the Riemmanian metric obtained from IG is the covariance matrix between the sufficient statistics\cite{Nielsen11,Caticha15}.
The EntDyn obtained in \cite{Pessoa20} is a diffusion process in the exponential family statistical manifold. 

Here we focus on a particular consequence of the EntDyn presented in \cite{Pessoa20}: the geometrical reciprocity in the linear regime. Meaning, as the dynamics evolves, the time derivative of the expected value of parameters is proportional to the gradient of entropy with reciprocal (symmetric) values. In the present article, we show how this reciprocity is a consequence of IG. 
This relates to what is known, in physics, as Onsager reciprocity principle (ORP) \cite{Onsager31} -- also referred to as Onsager reciprocal relations -- for non-equilibrium statistical mechanics. A brief review of ORP is presented in section 4 below.
Moreover, a diffusion process that reduces to a reciprocal relation in the linear regime is known in physics as the Onsager-Machlup process \cite{Onsager53}.  Unlike in EntDyn, the reciprocity found by Onsager is a consequence of the time reversal symmetry in an underlying dynamics, hence ORP is not trivially generalizable to dynamical processes beyond physics.

Gibbs distributions arise in equilibrium statistical mechanics as the maximum entropy distributions obtained when one chooses the conserved quantities in the Hamiltonian dynamics as sufficient statistics \cite{Jaynes57}. The generality of maximum entropy and the concept of sufficiency ensures that the exponential family is not limited to thermodynamics.
Since both the information theory and the IG approach guiding EntDyn are consequences of fundamental concepts in probability theory,  the reciprocity presented here is applicable for dynamical systems beyond physics. 
Just like ORP found applicability in non-equilibrium statistical mechanics\footnote{
It is also relevant to say that information theory ideas have been applied to nonequilibrum statistical mechanics (see \cite{Jaynes79} and posterior literature on maximum caliber e.g. \cite{Presse13})  which, unlike EntDyn, are approaches based on Hamiltonian dynamics. }, awarding Onsager with a Nobel prize in 1968, it is remarkable that a similar principle arises from a geometric framework.

The layout of the present article is as follows:
In the following section we will establish notation and give a brief review on exponential family/Gibbs distributions and their geometric structure endowed by IG.
In section 3 we will review the assumptions and results in \cite{Pessoa20}, obtaining the transition probabilities and expressing  EntDyn as a diffusion process on the statistical manifold.
In section 4 we comment on the similarities and differences in the reciprocity emerging from EntDyn and the one found by Onsager.

\section{Background - Exponential family}
We start by defining the exponential family of distributions $\rho$ for a given set $\mathcal{X}$ with elements $x \in \mathcal{X}$:
\begin{equation}
\rho(x|\lambda) \doteq \frac{q(x)}{Z(\lambda)}\exp[- \sum_{i=1}^n \lambda_{i}a^{i}(x)] \ ,\label{canonicaldefinition}%
\end{equation}
which can be obtained either as a result of maximization of entropy (as in \cite{Pessoa20}) or as a consequence of sufficiency (as in \cite{Nielsen11}). In \eqref{canonicaldefinition} the functions $a^i$ -- indexed up to a finite number $n$ -- are the sufficient statistics, the set of real values $\lambda = \{ \lambda_i \}$ are the parameters for the exponential family (or the Lagrange multipliers arising from maximization of entropy), $q(x)$ is an underlying measure (or the maximum entropy prior) and $Z(\lambda)$ is a normalization factor computed as
\begin{equation}\label{Zdefinition}%
Z(\lambda)=\int \dd x \ q(x)\exp[-  \lambda_{i}a^{i}(x)] \ .
\end{equation}
Here and for the remainder of this article we use the Einstein's summation notation, $A^i B_i = \sum_{i=1}^n A^i B_i $.
Details on how many well-known distributions can be put in the exponential form \eqref{canonicaldefinition} are found in \cite{Nielsen11}.

The expected values of sufficient statistics $a^i(x)$ in the exponential family can be found using 
\begin{equation}
    A^i \doteq \langle a^i(x)\rangle = -\pdv{}{\lambda_i} \log Z(\lambda) \ .
\end{equation}
Moreover, the entropy for $\rho$ in \eqref{canonicaldefinition} with a prior given by $q(x)$ is
\begin{equation} \begin{split} \label{GibbsEntropy}
     S(A) \doteq S[\rho|q] &= - \int dx \ \rho(x|\lambda(A)) \log\frac{ \rho(x|\lambda(A))}{q(x)} \\ &=  \lambda_i(A) A^i +\log Z (\lambda(A)) \ .
\end{split}\end{equation}
Thus, the entropy of the exponential family, as a function of the expected values $A = \{ A^i\}$, is the Legendre transform of $-\log Z$. It also follows that 
\begin{equation}\label{gradient}
    \lambda_i = \pdv{S}{A^i} \ .
\end{equation}

Each set of expected values $A$ corresponds to a single $\lambda$ per \eqref{gradient} and, because of this, we will use the expected values $A$ as coordinates for the space of probability distributions $\rho(x|\lambda(A))$ written, for simplicity, as just $\rho(x|A)$. IG consists of assigning a Riemmanian geometry structure to the space of probability distributions. Meaning, if we have a family of probability distributions parametrized by a set of coordinates, $P(x|\theta)$ where, $\theta = \{\theta^i\}$, the distance $\dd \ell$ between the neighbouring distributions $P(x|\theta+\dd \theta)$ and $P(x|\theta)$ is given by $\dd \ell^2 = g_{ij} \dd \theta^i \dd \theta^j$ where $g_{ij}$ is the Riemmanian metric. For probability distributions the only metric that is appropriate\footnote{The formal argument here is that, as proven by Cencov \cite{Cencov81}, FRIM is the only metric that is invariant under Markov embeddings, \change{ see also \cite{Ay17}}. Therefore, it is the only metric leading to a geometry consistent with the grouping property of probability distributions.} 
is the Fisher-Rao information metric (FRIM)

\begin{equation}
\label{firstmetric}
    g_{ij} = \int \dd x \ P(x|\theta) \frac{\partial\log P(x|\theta)}{\partial \theta^i}\frac{\partial\log P(x|\theta)}{\partial \theta^j}  .
\end{equation}
For the exponential family -- $\theta = A$ and $P(x|\theta) = \rho(x|A)$ in \eqref{firstmetric} --  FRIM can be expressed in the useful form 
\begin{equation}\label{HessianMetric}
    g_{ij} =- \pdv{S}{A^i}{A^j} \ .
\end{equation}
Having the exponential family and its IG established, we will explain the EntDyn obtained from it.

\section{Entropic dynamics on the exponential family manifold}

The dynamical process in EntDyn consists of evolving the probability distribution for a system in an instant $\rho(x|A)$ to the distribution describing a later instant $\rho(x'|A')$. 
This is done by assigning a probability $P(A)$ from which the system's state probabilities can be obtained as
\begin{equation}
    P(x) = \int \dd A \ P(A) \rho(x|A) \ ,
\end{equation}
which means that the conditional probability of $x$ given $A$ has to be of the exponential form defined in \eqref{canonicaldefinition}. Under EntDyn, we find the transition probabilities $P(A'|A)$ through the methods of information theory -- maximizing an entropy. 
As we will see below, the constraint presented establishes that the motion does not leave the manifold of exponential  distributions.

\subsection{Obtaining the transition probability}

\subs{The entropy} we opt to maximize must be so that the dynamical process accounts for both the change in state of the system -- evolving from $x$ to $x'$ -- and the change in the underlying distribution -- from  $\rho(x|A)$ to $\rho(x'|A')$. Hence,  it is the entropy for the joint transition, $P(x',A'|x,A)$, given by
\begin{equation}
\mathsf{S}[P| Q] = - \int \dd x'  \dd A' \ P(x', A' | x, A) \log(\frac{P(x', A' | x, A)}{Q(x', A' | x, A)}) \ ,
\label{DynamicalS}
\end{equation}
where $Q$ is the prior which we will determine below. Since the maximization of $\mathsf{S}$ will define the dynamics, we will  call it the dynamical entropy. 
One should not confuse the entropy for the exponential distribution, $S$ in \eqref{GibbsEntropy}, with the dynamical entropy, $\mathsf{S}$ in \eqref{DynamicalS}.

\subs{The prior} for the maximization of \eqref{DynamicalS} must implement a  continuous motion -- {meaning implementing short steps so that $A'$ converges in probability to $A$}. As explained in \cite{Pessoa20}, the least informative prior that does so is 
\begin{equation}
Q(x', A' |x, A)  \propto  g^{1/2}(A') \exp(-\frac1{2\tau}  g_{ij} \Delta A^i \Delta A^j) \ ,
\label{eq:EDGPrior}
\end{equation}
where $\Delta A^i = {A'}^i-A^i$ and $g=\det g_{ij}$. The parameter $\tau$ establishes continuity -- $\Delta A$ converges to 0 in probability when $\tau \rightarrow 0$. Later we will see how $\tau$ takes the role of time. 

\subs{The constraint} for the maximization of \eqref{DynamicalS} that implements a motion that does not leave the exponential manifold  is 
\begin{equation}
P(x', A' | x, A) = P(x' | x, A, A') P(A' | x, A) = \rho(x' | A') P(A' | x, A) \ . 
\label{eq:EDGConst}
\end{equation}
which means that the posterior distribution for $x'$ has to be conditioned only on $A'$ and to be in the exponential form \eqref{canonicaldefinition}.

\subs{The transition} that maximizes \eqref{DynamicalS} with the prior \eqref{eq:EDGPrior} and under \eqref{eq:EDGConst} is

\begin{equation} \label{transition1}
    P(x',A' | x,A) \propto \rho(x|A') e^{-S(A)} g^{1/2}(A') \exp( -\frac1{2\tau}  g_{ij} \Delta A^i \Delta A^j) \ .
\end{equation}
In order to obtain the transition probability from $A$ to $A'$, we must identify that from conditional probabilities in \eqref{eq:EDGConst} we have $P(A'|x,A) = \frac{P(x,A'|x,A)}{\rho(x'|A')}$ then apply a linear approximation of $S(A')$ near $S(A)$. As explained in \cite{Pessoa20} this results in
\begin{equation} \label{transition2}
P(A' | A) = \frac{g^{1/2}(A')}{\mathcal{Z}(A)} \  \exp(\pdv{S}{A^i} \Delta A^i -\frac1{2\tau} g_{ij} \Delta A^i \Delta A^j) \ ,
\end{equation}
where $\mathcal{Z}(A)$ is a normalization factor.
In the following subsection, we will explain how this transition probability results in a diffusion process on the statistical manifold.

\subsection{Entropic dynamics as a diffusion process}

Since according to \eqref{eq:EDGPrior} $\tau \rightarrow 0$ for short steps, we will compute the moments for $\Delta A$ in \eqref{transition2} up to order $\tau$ obtaining, as in \cite{Pessoa20},
\begin{equation}\label{moment}
            \expval{\Delta A^i} = \tau  g^{ij}\pdv{S}{A^j} - \frac{\tau}{2}\Gamma^i  + o(\tau) \ ,  
\end{equation}
where the upper indexes in $g$ to represent the inverse matrix, $ g_{ij} g^{jk}  = \delta_i^k$,  $\Gamma^i =  \Gamma^i_{jk} g^{jk} $ and $\Gamma^i_{jk}$ are the Christoffel symbols. In \cite{Pessoa20} we also show that the second and third moments of $\Delta A$ are
\begin{equation} \label{moment23}
        \expval{\Delta A^i\Delta A^j} = \tau g^{ij} + o(\tau) \   \quad \text{and} \quad \ \expval{\Delta A^i \Delta A^j \Delta A^k} =  o(\tau)  \ .
\end{equation}

The rules for a change from $A$ to $A'$ must be the same for a posterior change from $A'$ to $A''$. To keep track of the accumulation of changes, we will design the transitions as a Markovian process.
Time is introduced so that it parametrizes different instants, meaning $P(A) = P_t(A)$ and $P(A') = P_{t'}(A)$ for $t'>t$. Two important consequences of this dynamical design are: (i) time is defined as an emerging parameter for the motion -- the system is its own clock -- and (ii) the dynamic is Markovian by design, rather than proven to be a Markov chain in an existing time.
The way to assign the time duration, $\Delta t = t'-t$, so that motion looks the simplest is in terms of its fluctuations, meaning
$\Delta t \doteq \tau \propto  g_{ij} \Delta A^i \Delta A^j$.

Using $\tau$ as the time duration, the moments calculated in \eqref{moment} and \eqref{moment23} are what is known as smooth diffusion (see e.g. \cite{Nelson85}). Therefore the dynamics can be written as a Fokker-Planck (diffusion) equation which written in invariant form is (see \cite{Pessoa20})
\begin{equation}\label{FPeq}
\pdv{p}{t} = - \frac1{g^{1/2}}  \pdv{A^i}\qty(g^{1/2} p  v^i)\ ,
\qq{where}
v^i =  g^{ij} \pdv{S}{A^j} -   \frac{g^{ij}}{2p}\pdv{p}{A^j}  \ ,
\end{equation}
and $p$ is the invariant probability density $p(A) \doteq \frac{P(A)}{\sqrt{g(A)}}$.
The first term for $v^i$ in \eqref{FPeq} is the diffusion drift velocity -- related to the average motion in the diffusion process -- while the second term is the osmotic velocity -- guided by the differences in probability density.
The drift velocity is linear with respect to the gradient of entropy, $\pdv{S}{A^i}$, with symmetric proportionality factors -- since the metric tensor is symmetric -- hence, we say that the drift velocity is reciprocal to the gradient of entropy. 
In the following section we will make further comments on  reciprocity in both EntDyn and ORP.

\section{Reciprocity}
In Onsager's approach to nonequilibrim statistical mechanics \cite{Onsager31}, we suppose that a thermodynamical system is fully described by a set of parameters $\xi$ evolving in time $\xi(t) = \{ \xi^i(t)\}$. ORP states that, near an equilibrium point $\xi_0$, the rate of change for $\xi$ is reciprocal to the entropy gradient, meaning
\begin{equation} \label{ORPeq}
    \dv{\xi^i}{t} =  \gamma^{ij} \pdv{S}{\xi^j}
\end{equation}
where $\gamma^{ij}=\gamma^{ji}$ and $S$ is the thermodynamical entropy -- \change{in mathematical terms it is} the Gibbs distribution entropy \eqref{GibbsEntropy} written in terms of $\xi$.

In ORP, the equilibrium value, $\xi_0$, has to both be a local maxima of entropy and a fixed point for the dynamics of $\xi$. That means
\begin{equation}
    \dv{\xi^i}{t} \at{\xi=\xi_0} = 0 \ \ \ \text{and} \ \ \ \pdv{S}{\xi^i} \at{\xi=\xi_0} = 0 \ .
\end{equation}
If we make a linear approximation for both $S$ and $\dv{\xi^i}{t}$ around $\xi_0$ we obtain
\begin{equation}\label{eq19}
    \dv{\xi^i}{t} = - L^i_k \beta^{kj}  \pdv{S}{\xi^j} \ ,
\end{equation}
where 

\begin{equation} \label{Landbeta}
    L^i_j \doteq \pdv{}{\xi^j}\dv{\xi^i}{t}\at{\xi=\xi_0}  \qq{and} \beta_{ij} \doteq - \frac{\partial^2 S}{\partial \xi^i \partial \xi^j} \at{\xi=\xi_0} \ .
\end{equation}
which, when compared to \eqref{ORPeq}, yields $\gamma^{ij}  = - L^i_k \beta^{kj} $.
Hence proving ORP equates to proving $ L^i_j = L^j_i $, since $\beta$ is already symmetric\footnote{ Note that although in \eqref{HessianMetric} $\beta$ matches FRIM  for $\xi = A$, one cannot say that $\beta$ is FRIM in general since \eqref{HessianMetric} does not transform covariantly for an arbitrary change of coordinates $A \rightarrow \xi(A)$.} for being defined as a second derivative \eqref{Landbeta}. 
The usual proof of the symmetry of $L$ (see e.g.  \cite{Landau80}) consists of assuming that the fluctuations of $\xi$ are symmetric through time reversal, meaning
$\expval{ \xi^i(0) \xi^j(t)} = \expval{\xi^i(t) \xi^j(0) }$. Unlike the formalism derived by EntDyn, this cannot be expected to be general for dynamical systems.

A direct comparison with \eqref{moment} means that when the system is described by the expected value parameters, $\xi =A$, EntDyn reduces, in the first moment, to a motion that is reciprocal to the gradient of entropy. The  second term for $\expval{\Delta A^i}$ in \eqref{moment} is a correction for motion that probes the manifold's curvature. Also, as discussed, in \eqref{FPeq} we see that the diffusion process resulting from EntDyn has drift velocity reciprocal to the gradient of entropy, analogous to the Onsager-Machlup process \cite{Onsager53}.

Although it would make a trivial connection between EntDyn and nonequilibrium statistical mechanics, one might find the rates of change described by it to be a physically unrealistic model since substituting $\gamma^{ij}$ by $ g^{ij}$ for $\xi = A$ in \eqref{ORPeq} and comparing to \eqref{eq19} implies:
\begin{equation} \label{motionthatfitsxi}
L^i_j = - \ \delta^i_j\
\quad  \rightarrow \quad
\dv{\xi^i}{t} \propto  \ (\xi_o^i - \xi^i) \ .
\end{equation}
The right hand side above is obtained by substituting the left hand side into the definition of $L$ in \eqref{Landbeta} and integration. 
In thermodynamics,  one can not expect macroscopic parameters, such as internal energy and number of particles, to evolve independently. 

\section{Conclusion and perspectives}

An entropic framework for dynamics on the exponential family manifold was presented here. The change in the expected values in \eqref{moment} is reciprocal to the gradient of entropy in the manifold, while \eqref{FPeq} gives a drift velocity that is also reciprocal to the gradient of entropy. The proportionality factors in this reciprocal motion is given by the inverse FRIM, therefore reciprocity is a consequence of IG.

From \eqref{motionthatfitsxi}, one should not consider the dynamics developed here as a method for nonequilibrum statistical mechanics. The reciprocity found here does not necessarily substitute ORP. 
That said, EntDyn offers a systematic method to find dynamics aligned with information theory, while Onsager’s approach is based solely on calculus considerations around a supposed fixed point. Another way to present it is to say that ORP is based on an understanding of thermodynamics guided by Hamiltonian dynamics, while the EntDyn developed here is inspired by the generality of exponential families, not relying on an underlying Hamiltonian dynamics. Both arrive at reciprocal relations. 


\vspace{.75cm}
\subs{Acknowledgments}
I thank N. Caticha, A. Caticha, F. Xavier Costa, B. Arderucio Costa, and M. Rikard for insightful discussions in the development of this article. 
I was financed in part by  CNPq -- Conselho Nacional de Desenvolvimento Científico e Tecnológico-- (scholarship GDE 249934/2013-2).


\bibliographystyle{splncs03-num}
\bibliography{referencias}

\end{document}